\documentclass[a4paper,11pt]{article}
\usepackage[utf8]{inputenc}

\usepackage{amsmath}
\usepackage{amsfonts}
\usepackage{amssymb}

\usepackage[T1]{fontenc}
\usepackage{gensymb}
\usepackage{lmodern}
\usepackage{mdframed}

\usepackage[hyphens]{url}
\usepackage{hyperref}
\usepackage[hyphenbreaks]{breakurl}

\usepackage[left=3cm,right=3cm,top=3cm,bottom=3.2cm]{geometry}

\date{1 September 2020}
\author{
Maurice Chiodo\footnote{King's College, University of Cambridge. \texttt{mcc56@cam.ac.uk} .} \ \
  Dennis M\"uller\footnote{RWTH Aachen University \texttt{dennis.mueller3@rwth-aachen.de} .}   
}

\title{Questions of Responsibility: Modelling in the Age of COVID-19}
\begin{document}

\maketitle

\let\thefootnote\relax\footnotetext{2020 \textit{AMS Classification:} 01A80. 00A30.}
\let\thefootnote\relax
\let\footnote\relax\footnotetext{Appeared in \textit{SIAM News 53, No. 7, 1 September 2020}. Available online at 
\url{https://sinews.siam.org/Details-Page/questions-of-responsibility-modelling-in-the-age-of-covid-19}.}

\section*{}
Throughout history, humanity has always battled pandemics with tools that are available in the time period. Today, these tools involve a combination of big data techniques and traditional epidemiological models. In the ongoing fight against COVID-19, many mathematicians work at the \textit{third} level of ethical engagement: \textit{taking a seat at the table of power} (we introduced \textit{SIAM News} readers to the four levels of ethical engagement for mathematicians in 2018) \cite{Chiodo.2018}. When conducting COVID-19-related research, mathematicians may find themselves in public service roles. They might engage in policy discussions, collaborate with politicians and other non-scientists, and use the predictive power of mathematics to shape substantial decisions. We aim to promote good practice for ethical mathematical modelling in these difficult times. We begin with some observations about mathematical modellers, then derive eight questions that all mathematicians should constantly ask themselves.

During a pandemic, circumstances demand rapid and immediate results. Consequently, researchers often publish or communicate preliminary findings before they are properly peer-reviewed. Like all scientists, mathematicians must strike a careful balance between quality and speed. Big data, analytics, and mathematical modelling have become the go-to options for many researchers in the age of quantitative decision-making. In particular, some may think that mathematical quantification \textit{always} improves a situation. However, believing that “any model is better than no model” is potentially dangerous when the results influence life-or-death decision-making in areas ranging from high-level political offices to hospital frontlines.

Alternative methods of problem-solving do exist, and providing a false sense of security can skew other people’s decision spaces. Scientists and mathematicians carry a burden of responsibility when they apply their tools and models to better understand a worldwide pandemic. Politicians like to be “led by science,” but who is held responsible for the ensuing decisions? Is it the mathematicians, the politicians, or both? Generally, mathematicians are not politicians and politicians are not mathematicians. In order to work and communicate effectively as a team, each group might need to learn a bit about the other side. And if mathematicians are engaging with the political decision process, what exactly should they be doing?

Proper science communication is a crucial component of a mathematician’s role, now more than ever. For example, the British government recently improperly explained its new alert level system to the general public as “COVID alert level = rate of infection + number of infections” \cite{Woolley.13092022}. This statement both confused people and lacked actual information, a combination that can seriously erode the public’s trust in mathematics. When engaging in and supporting higher-level decision-making, mathematicians must find adequate ways to communicate their results and decisions to the public. Some researchers are acting with the utmost transparency, including a group of self-organized scientists called the Independent SAGE$^1$\let\thefootnote\relax\footnotetext{$^1$\url{https://www.independentsage.org/}} that aims to provide transparent government advice during COVID-19.

Mathematical models, which range from predictions for reopening schools and universities to statistics pertaining to children’s viral loads, can incur real-life costs. Therefore, clear and consistent advice is necessary for leaders to make proper decisions. For instance, a school principal cannot prepare his/her school for reopening when predictions—and thus government advice—are constantly changing or contradictory \cite{Pidd.06122020}. In all of their equations and theorems, mathematicians cannot forget the people at the frontline and those making the political decisions. They must ensure that their work is practical, useful, and does not overly constrain politicians’ decision spaces. Otherwise, a lack of usability—combined with public distrust and an anti-science bias—could erode mathematics’ authority \cite{Howard.18.6.2020}.

The COVID-19 pandemic has also underlined mathematicians’ increasing levels of political responsibility, which extend beyond their technical work. Private actions can undermine important mathematical and scientific research, even when the studies are sound. Such was the case when epidemiologist and mathematician Neil Ferguson resigned from his U.K. government advisory position after violating lockdown rules \cite{BBC.05052020}.

While assisting other mathematicians on COVID-19 modelling efforts, we observed that splitting large problems into smaller, manageable subproblems and subsequently compartmentalising teams consistently led to communication issues. Teams need members who are dedicated to ensuring seamless communication and avoiding structural problems in the models. Furthermore, identifying only one member to address ethical issues within a large team is insufficient, as there is simply too much ground to cover. Every mathematician must, to some extent, be aware of his/her own moral responsibilities.

Mathematicians never work in a political or ethical vacuum, and this is particularly true when their efforts impact a pandemic response. In these times, they cannot simply expect their work to “speak for itself.” Instead, mathematicians must vocalise the importance and relevance of their research and ensure that they possess a solid understanding of their own ethical standards and the consequences of their work. For example, data scientist Rebekah Jones was reportedly fired from her job with Florida’s Department of Health after refusing to manipulate COVID-19 data on a public-facing portal because she considered the requests unethical \cite{Iati.13.6.2020}. She has since begun publishing data on her private website.

Understanding, processing, and inferring meaning from complex data requires a broad perspective. An emphasis on the reproducibility of modelling results, which enables “outsiders” with different views to scrutinize the work in question, is also essential. When software engineers identified sloppy code in a British COVID-19 model, journalists were quick to report on it. However, it is now evident that these results are actually reproducible \cite{SinghChawla.2020}.

But shouldn't mathematicians aim for more than just clean code? Computer scientists have long understood the value of commenting source code, but how common is it for mathematicians to properly comment and document their models? For critical results, mathematicians must document the \textit{reasoning} that inspired their solution. Over the last decade, computer science and other disciplines have developed tools to better document data and models for artificial intelligence (AI) production systems. With slight adaptations, existing tools like Datasheets for Datasets \cite{Gebru.24032018} or Model Cards for Model Reporting \cite{Mitchell.2019} can help mathematicians make their models and decisions easier for the public to digest and scrutinise. These tools may also help ensure that researchers are not working with biased or incomplete datasets.

Good practice in AI development expects scientists to consider the origins of and biases in their datasets; document their models, data, and impact; consider feedback loops; and remain aware of their responsibilities to the public. These principles also apply in pandemic modelling. Mathematicians must not fall victim to limitations or biases in their datasets or thinking. As a rough guideline for others, we summarise our observations in eight questions that mathematicians who are working on COVID-19 models should keep in mind:

\begin{enumerate}
\item Am I using \textit{authorised and morally obtained datasets}?
\item Do my co-workers, superiors, and I have \textit{sufficient perspective}, and do we understand the \textit{limitations and biases in our data and thinking}?
\item What \textit{optimisation objectives and constraints} have I chosen, and what are their \textit{real-life costs}?
\item Am I properly considering how to \textit{comment and document my model and communicate the results} to those who need them?
\item Are my results \textit{explainable and falsifiable}?
\item Do I have techniques to handle \textit{feedback loops} and the \textit{large-scale impact} of my work?
\item Am I aware of other \textit{non-mathematical aspects} and the \textit{political nature} of my work?
\item Do I have a \textit{non-technical response strategy for when things go wrong}? Do I have \textit{peers who support me with whom I can talk}?
\end{enumerate}

Of course, these aspects are not entirely new. Emanuel Derman and Paul Wilmott laid out their “Financial Modelers’ Manifesto” in the wake of the 2008 global financial crash \cite{Dermann.2009}. But mathematicians should not overlook the tendency to forget basic principles in the heat of the moment. We designed these eight questions to promote good practices for ethical modelling. Good practices value everyone — colleagues and the model subjects alike. During these arduous times, mathematicians become part of bigger processes when they use their expertise to help policymakers and frontline workers make decisions. Contributing to the COVID-19 response is undeniably challenging; the stakes are high, the situation is unfamiliar, information is uncertain, and circumstances change quickly. Mathematicians must formulate their work responsibly and look beyond their LaTeX files to ensure that their work is adequately communicated and well-received among scientists, decision-makers, and the general populace.

\par\noindent\rule{\textwidth}{0.4pt}
\textit{For the past four years, the Cambridge University Ethics in Mathematics Project has been examining the ethical issues and social ramifications of all forms of mathematical research. Our work—which includes articles, two organised conferences, and a list of talks—is available online.}$^2$\let\thefootnote\relax\footnotetext{$^2$\url{https://www.ethics.maths.cam.ac.uk/}}\textit{ We hope that our insights and tools prove useful to mathematicians working on COVID-19, and welcome any contact or requests for collaboration, assistance, or advice. The mathematical community works better when it works together.}


\end{document}